\documentclass[runningheads,a4paper]{llncs}

\usepackage{amssymb}
\usepackage{amsmath}
\usepackage{graphicx}

\usepackage{url}
\newcommand{\keywords}[1]{\par\addvspace\baselineskip
\noindent\keywordname\enspace\ignorespaces#1}

\newcommand{\Gal}{\mathrm{Gal}}
\newcommand{\fpbar}{\overline{\bbbf}_p}
\newcommand{\Frob}{\mathrm{Frob}}
\newcommand{\GL}{\mathrm{GL}}

\begin{document}

\mainmatter  

\title{Enumerating Galois Representations in Sage}


%
%
\author{Craig Citro\inst{1}\inst{2} \and Alexandru Ghitza\inst{3}}

%

\institute{Google, Seattle WA\\
\and
Department of Mathematics, University of Washington,\\
Box 354350, Seattle WA 98195-4350, USA,\\
\email{craigcitro@gmail.com}
\and
Department of Mathematics and Statistics, University of Melbourne,\\
Parkville VIC 3010, Australia,\\
\email{aghitza@alum.mit.edu}}

%
%

\maketitle

\begin{abstract}
  We present an algorithm for enumerating all odd semisimple
  two-dimensional mod $p$ Galois representations unramified outside $p$.  
  We also discuss the implementation of this algorithm in Sage and give 
  a summary of the results we obtained\footnote{The authors wish to
  thank Kevin Buzzard for providing several corrections and a significant
  improvement to Theorem~\ref{thm:comparison}, and the referees for 
  suggesting improvements to the exposition.}.
\end{abstract}

\section{Introduction}

A great deal of arithmetic questions have found natural interpretations
(and often, answers) within the realm of Galois representations and 
modular forms: such applications include Diophantine equations,
quadratic forms, or the study of combinatorial-arithmetic objects such
as partitions.  In this context, it is of interest to dispose of computational
tools for working with modular forms and Galois representations.

In this note, we focus on two-dimensional \emph{Galois representations mod~$p$}, 
i.e. continuous group homomorphisms
\begin{equation*}
  \rho\colon\Gal(\overline{\bbbq}/\bbbq)\longrightarrow\GL_2(\fpbar)\enspace.
\end{equation*}
(More precisely, we consider such representations which are semisimple,
unramified outside $p$, and odd.  For the theoretical background, we refer the
reader to Khare's survey~\cite{Khare2007} or to Edixhoven's 
paper~\cite{Edixhoven1992}.)

By Serre's conjecture, now a theorem of Khare-Wintenberger 
(see~\cite{Khare2009a},~\cite{Khare2009b}), these representations
are closely related to modular forms (mod $p$) of level $1$ which
are eigenvectors for all the Hecke operators.
If $f$ is such a form, of weight~$k$ and 
eigenvalues $(a_\ell)$, then for all primes $\ell\neq p$ we have
\begin{equation*}
  \mbox{charpoly}(\rho(\Frob_\ell))=X^2-a_\ell X+\ell^{k-1}\enspace,
\end{equation*}
where $\Frob_\ell$ is a Frobenius element at $\ell$ inside 
$\Gal(\overline{\bbbq}/\bbbq)$.

The \emph{(Hecke) eigensystem} corresponding to a
mod $p$ eigenform $f$ is the sequence $\left(a_\ell\right)$ 
of eigenvalues indexed by all primes~$\ell\neq p$.  The $i$-th \emph{twist} of 
$\left(a_\ell\right)$ is by definition the
eigensystem $\left(\ell^i a_\ell\right)$.  We write~$[a_\ell]$ for a finite
truncation of $(a_\ell)$, where the cutoff point will be clear from the
context.

Inspired by a remark of Khare\footnote{From Sect.~8 of~\cite{Khare2007}:
``[\dots] there are only finitely many semisimple $2$-dimensional mod~$p$ 
representations of $\Gal(\overline{\bbbq}/\bbbq)$ of bounded (prime-to-$p$
Artin) conductor. \emph{It will be of interest to get quantitative refinements 
of this}.''}, we have set out to enumerate
all odd semisimple mod $p$ representations which are unramified outside 
$p$.  This corresponds to enumerating all the Hecke eigensystems
which occur in spaces of level~$1$ modular forms mod $p$.

\section{Description of the Algorithm}
The starting point is a classical result in the theory of modular forms
mod $p$ (see Theorem~3.4 in~\cite{Edixhoven1992}): every Hecke eigensystem 
occurs, up to twist, in weights less than or equal to~$p+1$.  Therefore it 
suffices to generate the spaces
$M_k$ for weights $4\leq k\leq p+1$ and find all the eigenforms in them,
which will produce all the Hecke eigensystems up to twist.  This list
may however contain duplicates;
we investigate this question in detail in~\cite{CitroGhitza},
where we prove

\begin{theorem}\label{thm:comparison}
  \
  \renewcommand{\labelenumi}{(\alph{enumi})}
  \begin{enumerate}
    \item Let $f_1$ and $f_2$ be eigenforms of weights $k_1,k_2\leq p+1$.
      If $f_1$ and $f_2$ have the same eigensystem up to twist, then 
      $k_1+k_2=p+1$ or $k_1+k_2=p+3$.
    \item Let $f_1$ and $f_2$ be eigenforms of weights related in one of
      the ways described in (a).  If $f_1$ and $f_2$ do not have the same
      eigensystem up to twist, then this is detected by a prime 
      $\ell\neq p$ satisfying $\ell\leq (p+1)/6$.
  \end{enumerate}
\end{theorem}

In the process of proving Theorem~\ref{thm:comparison}, we obtained
the following lower bound, which improves the best known lower bound
(due to Serre, see Sect.~8 in~\cite{Khare2007}) by a factor of two:

\begin{theorem}
  Let $p>19$ be prime.
  The number of odd semisimple $2$-dimensional Galois representations
  mod $p$ which are unramified outside $p$ is bounded below by $p(p-1)/2$.
\end{theorem}

\subsection*{Algorithm: Enumerate Galois representations mod $p$
up to twist}

\begin{enumerate}
  \item For $4\leq k\leq p+1$:
    \begin{enumerate}
      \item Compute a basis for the space $M_k$.
      \item Decompose the space into Hecke eigenspaces.
      \item For each eigenform, compute the eigenvalue $a_\ell$ of 
        $T_\ell$ for primes $\ell$ up to the bound from 
        Theorem~\ref{thm:comparison}.  Store 
        $\left(k,\left[a_\ell\right]\right)$.
      \end{enumerate}
  \item Remove duplicates: given $\left(k_1,\left[a_\ell\right]\right)$
    and $\left(k_2,\left[b_\ell\right]\right)$ such that $k_1+k_2=p+1$
    or $p+3$, check whether $\left[b_\ell\right]$ is a twist of 
    $\left[a_\ell\right]$.
\end{enumerate}

This creates the list of equivalence classes (up to twist) of Hecke
eigensystems mod $p$.  It is now straightforward to apply the twist 
operation to each list element and generate the list of all Hecke 
eigensystems.

\section{Sage Implementation and Results}

Our task requires computing the action of Hecke operators on spaces of
modular forms of high weight.  Sage~\cite{Sage} offers several 
implementations of these spaces for arbitrary levels.
We have initially
used modular symbols over finite fields for generating the
lists of eigenforms, but this method becomes quite slow as the weight
increases.  Restricting to level~$1$ allows us to take advantage 
of a much faster way of working with these spaces: the 
Victor Miller basis (see Sect.~2.3 in~\cite{Stein2007} for the properties
and the algorithm Sage uses to compute this basis).

We then use one Hecke operator $T_\ell$ at a time to decompose the
space $M_k$ into eigenspaces.  This requires (at most) the first $k/12$ 
primes $\ell$ (see the appendix of~\cite{Lario2002}).

We have run the Sage implementation of our algorithm for all primes up 
to~$211$ (see Table~\ref{table:results}). 
Apart from keeping track of the number of equivalence classes of
eigensystems and the total number of eigensystems, we save the list of
equivalence classes; given this it is very easy to take twists and
generate the entire list.  

\vspace{-0.5cm}

\begin{table}[h]
  \centering
  \begin{tabular*}{0.95\textwidth}{@{\extracolsep{\fill}}r|r||r|r||r|r||r|r||r|r||r|r}
    $p$ & number & $p$ & number & $p$ & number & 
    $p$ & number & $p$ & number & $p$ & number\\ \hline 
    $2$ & $1$ &
    $23$& $264$ & 
    $59$& $4234$ &
    $97$& $19200$ &
    $137$& $53992$ &
    $179$& $119705$ \\
    $3$ & $1$ & 
    $29$& $532$ & 
    $61$& $4800$ & 
    $101$& $21600$ & 
    $139$& $55752$ &
    $181$& $124020$ \\
    $5$ & $4$ & 
    $31$& $630$ & 
    $67$& $6237$ & 
    $103$& $22797$ &
    $149$& $69264$ &
    $191$& $145445$\\
    $7$ & $9$ & 
    $37$& $1044$ & 
    $71$& $7420$ & 
    $107$& $25546$ &
    $151$& $71700$ &
    $193$& $150144$ \\
    $11$& $35$ & 
    $41$& $1480$ & 
    $73$& $8136$ & 
    $109$& $27216$ & 
    $157$& $80340$ &
    $197$& $160132$ \\
    $13$& $48$ & 
    $43$& $1701$ & 
    $79$& $10257$ &
    $113$& $30240$ &  
    $163$& $90477$ &
    $199$& $164637$ \\
    $17$& $112$ &
    $47$& $2185$ & 
    $83$& $12054$ &
    $127$& $42903$ & 
    $167$& $97276$ &
    $211$& $196560$ \\
    $19$& $153$ &
    $53$& $3172$ &
    $89$& $14784$ &
    $131$& $46735$ &
    $173$& $108016$\\
  \end{tabular*}
  \vspace{0.1cm}
  \caption{Number of Galois representations mod $p$}
  \label{table:results}
\end{table}

\vspace{-1.0cm}

Khare guesses in~\cite{Khare2007} that the number of Galois representations of 
the type we are considering should be asymptotic to $p^3/48$.  There are two 
phenomena that can contribute to the actual number being smaller than the guess: 
(i) the existence of ``companion forms'', which in our context appear as duplicate 
equivalence classes of eigenforms; (ii) the failure of ``multiplicity one'' for 
Hecke eigenvalues mod $p$, which results in some spaces $M_k$ not contributing 
their dimension's worth of eigenforms.  In the range of our computations, the
actual number of representations stays very close to the best known upper 
bound\footnote{For instance, for $p=211$ the quotient between the actual
number ($196560$) and the upper bound ($196665$) is about $0.9995$.},
suggesting that the two phenomena are indeed quite rare.  We expect this
trend to be confirmed by further computations.

\bibliographystyle{splncs03}
\bibliography{icms}

\end{document}